%% file: extcond.tex
\documentclass{amsart}
\usepackage{amssymb}
\usepackage{amscd}
\usepackage{amsmath}
\usepackage{amsthm}
\usepackage[dvips]{graphicx}
\usepackage[dvipdfm,%
bookmarks=true,%
bookmarksnumbered=true,%
]{hyperref}

\newcommand{\Z}{{\mathbb Z}}

\theoremstyle{plain}
\newtheorem{Thm}{Theorem}[section]
\newtheorem{Lem}[Thm]{Lemma}
\newtheorem{Prop}[Thm]{Proposition}

\theoremstyle{definition}
\newtheorem{Def}[Thm]{Definition}
\newtheorem{Ex}[Thm]{Example}
\newtheorem{Rem}{Remark}[section]

\title
[External edge condition and group cohomologies]
{External edge condition and \\ group cohomologies 
associated with \\ the quantum Clebsch-Gordan condition}
\author{Hajime Fujita}
\date{}

\keywords{trivalent graph ; admissible coloring ; TQFT}

\subjclass[2000]{Primary 05C15; Secondary 58D19, 57R56}

\address{Department of Mathematics, Gakushuin University 1-5-1, 
Mejiro, Toshima-ku, Tokyo, 171-8588, Japan} 

\email{hajime@math.gakushuin.ac.jp}

\begin{document}
\maketitle

\begin{abstract}
In this article we determine the 
structure of a twisted first cohomology group of 
the first homology of a trivalent graph 
with a coefficient associated with the 
quantum Clebsch-Gordan condition. 
As an application we give a characterization of 
a combinatorial property, the external edge condition (\cite{hjfj}), 
which is defined in the study of the Heisenberg 
representation on the TQFT-module. 
\end{abstract}


\section{Introduction}

In \cite{hjfj} we computed representation matrices and  
characters of a Heisenberg representation 
in terms of {\it admissible colorings} of a given ribbon graph. 
The Heisenberg representation was constructed in the framework of 
(2+1)-dimensional topological quantum field theory 
by Blanchet {\it et al.} in \cite{BHMV} 
and Andersen {\it et al.} in \cite{AndMas} 
based on skein theoretical and algebro-geometrical settings respectively. 
In the computation of characters, a combinatorial property 
which we call the {\it external edge condition} plays 
an essential role. 
In fact as we noted in \cite[Remark~6.1]{hjfj}, 
for each twisted group 1-cocycle of the homology of the graph 
with a coefficient associated with the 
quantum Clebsch-Gordan condition 
we can define a Heisenberg representation, 
and if the cocycle satisfies the external edge condition, 
then the representation is isomorphic to that in \cite{BHMV} and \cite{AndMas}. In this article we show that the converse is true 
in a suitable functorial sense 
(Theorem~\ref{charofext}). 
One of the key observations 
is the description of the twisted first cohomology group 
in terms of stabilizers of the action of the homology of the graph 
on the set of admissible colorings 
(Theorem~\ref{twistcohom}).

This article is organized as follows. 
In Section~2 we prepare some notations and definitions including 
admissible colorings and the twisted cohomology group associated with 
the graph and admissible colorings. 
In Section~3 we determine the structure of the twisted first cohomology. 
In Section~4 we recall the definition of the external edge condition 
and show the existence of the external edge class. 
In Section~5 we give an alternative description of the twisted first cohomology 
in terms of the induced representation on the TQFT-module. 
Finally we introduce a category of graphs and 
give a characterization of the external edge class 
in the categorical setting. 

The author's research is partially supported by 
JSPS Grant-in-aid for Young Scientists (No.21840045).

\section{Notations and definitions}
In this article we fix a positive integer $k$ called the {\it level} 
and a commutative ring $R$ with unit. 
We denote by $R^{\times}$ the group of units. 

Let $\Gamma$ be a unitrivalent graph
\footnote{A unitrivalent graph is a graph 
whose vertices are trivalent or univalent. 
In this article we do not need any ribbon structure of the graph. 
See the comments before Theorem~\ref{existextedge} also.} 
with 
$3g-3+2n$ edges $\{f_l\}$, $2g-2+n$ 
trivalent vertices $\{v_i\}$ and 
$n$ univalent vertices $\{w_m\}$.  
We assume that first $n$ edges $f_1,\cdots,f_n$ 
have univalent vertices $w_1,\cdots,w_n$. 
Note that the first homology group $H_1(\Gamma)$ of $\Gamma$ with 
$\Z/2$-coefficient is isomorphic to $(\Z/2)^{g}$. 
Let $\vec{j'}=(j'_1,\cdots,j_n')
\in \{0,\frac{1}{2},\cdots,\frac{k}{2}\}^n$ be a 
pair of half integers. 
We call the pair $(\Gamma;\vec j')$ a 
colored graph. 

\begin{Def}[Admissible colorings]\label{QCG}
A labeling of the set of edges 
$\vec j=(j_l):\{f_l\}\to \{0,\frac{1}{2},\cdots,\frac{k}{2}\}$ is an 
{\it admissible coloring of level} $k$ for 
$(\Gamma;\vec{j'})$ 
if $\vec j(f_l)=j_l'$ for $l=1,\cdots n$ and 
the following condition, which is called 
the {\it quantum Clebsch-Gordan condition of level $k$}, 
is satisfied 
for each trivalent vertex $v_i$ 
with three edges $f_{i_1},f_{i_2}$ and $f_{i_3}$: 
$$
\left\{\begin{array}{lll}j_{i_1}+j_{i_2}+j_{i_3}\in \Z 
\\|j_{i_1}-j_{i_2}|\leq j_{i_3}\leq j_{i_1}+j_{i_2} 
\\ j_{i_1}+j_{i_2}+j_{i_3}\leq k. \end{array}\right. 
$$
Here we put $j_l:=\vec j(f_l)$. 
If a trivalent vertex 
$v_i$ has only two edges $f_{i_1}$ and $f_{i_2}=f_{i_3}$, 
then we interpret these conditions as the 
corresponding condition with $j_{i_2}=j_{i_3}$. 
Denote by $QCG_k(\Gamma;\vec{j'})$ 
the set of all admissible colorings of level $k$ 
for $(\Gamma;\vec{j'})$. 
\end{Def}

\begin{Def}
We define an action of $H_1(\Gamma)$ on 
$QCG_k(\Gamma;\vec{j'})$ by 
$$
\lambda:\vec j\mapsto\lambda\cdot \vec{j}
=\left(j_1,\cdots,\frac{k}{2}-j_l,\cdots,j_{3g-3}\right) \ 
({\rm  all} \  l \ {\rm with} \ f_l \ 
{\rm lying \ on} \ \lambda). 
$$
Denote by $\overline{QCG_k}(\Gamma;\vec{j'})$ 
the quotient set with respect to this action and 
by $[\vec j]$ the equivalence class of $\vec j\in QCG_k(\Gamma;\vec{j'})$. 
\end{Def}

\begin{Ex}
For the graph $\Gamma$ in Figure~\ref{H_1-action} 
and $\lambda=f_3+f_4+f_5+f_6\in H_1(\Gamma)$, the action of $\lambda$ on $QCG_k(\Gamma;\vec j')$ 
is given by; 
$$
\lambda:(j_3, j_4, j_5, j_6,j_7,j_8,j_9,j_{10})\mapsto 
\left(\frac{k}{2}-j_3, \frac{k}{2}-j_4, \frac{k}{2}-j_5, \frac{k}{2}-j_6, 
j_7,j_8,j_9,j_{10} \right). 
$$

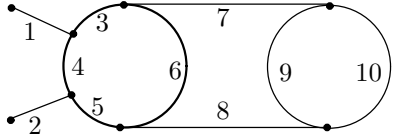
\begin{figure}[htb]
\begin{center}
\input{h1action.tex}
\caption{$\Gamma$ and $\lambda$.}\label{H_1-action}
\end{center}
\end{figure}
\end{Ex}

\begin{Def}
Put $A(\Gamma;\vec{j'}):=(R^{\times})^{QCG_k(\Gamma;\vec{j'})}$  
and define the action of $H_1(\Gamma)$ on 
$A(\Gamma;\vec{j'})$ 
by interchanging of entries, 
$\lambda:(c_{\vec{j}})\mapsto (c_{\lambda\cdot \vec{j}})$.
We denote by $C^*(\Gamma;\vec{j'})$ the twisted cochain group 
$C^*(H_1(\Gamma);A(\Gamma;\vec{j'}))$ 
and by $H^*(\Gamma;\vec{j'})$ its cohomology group 
$H^*(H_1(\Gamma);A(\Gamma;\vec{j'}))$. 
\end{Def}

\section{Structure of the first cohomology}

Our main interest is the first cohomology group
$H^1(\Gamma;\vec{j'})$. 
We first give the structure of it in terms of 
the homology of the given colored graph. 

\begin{Thm}\label{twistcohom}
For a colored  graph $(\Gamma;\vec j')$ 
we have the following isomorphism: 
$$
H^1(\Gamma;\vec{j'}) \ \cong
\bigoplus_{[\vec j]\in \overline{QCG_k}(\Gamma;\vec{j'})}
H_1(\Gamma)_{\vec{j}}^*, 
$$where we put 
$H_1(\Gamma)_{\vec{j}}:=\{\lambda\in H_1(\Gamma) \ | 
\lambda\cdot\vec j=\vec j\}$ and 
$(\cdot)^*:=Hom( \ \cdot \ ,R^{\times})$. 
\end{Thm}
To show Theorem~\ref{twistcohom} we show the following lemma.
\begin{Lem}\label{coboundary}
We have the following description of the 
first coboundaries $B^1(\Gamma;\vec{j'})$. 
$$
B^1(\Gamma;\vec{j'})=\{\delta=(\delta_{\vec{j}})\in Z^1(\Gamma;\vec{j'}) \ | \ 
\delta_{\vec{j}}(\lambda)=0 \ {\rm for \ all} \ \vec{j}\in 
QCG_k(\Gamma;\vec j') \ {\rm and} \ \lambda\in H_1(\Gamma)_{\vec j} \}.
$$
\end{Lem}
This lemma follows from the following general proposition.
\begin{Prop}
Let $G$ be an abelian group 
acting on a finite set $A$.  
Let $K$ be an abelian group and 
$\tau:G\to K^A$ ($\tau(g)=(\tau_a(g))$) a twisted 1-cocycle of $G$, 
where $K^A$ is the direct product and 
$G$ acts on $K^A$ by interchanging of entries, 
$g:(r_a)_{a\in A}\mapsto(r_{g\cdot a})_{a\in A}$ \ $(g\in G$). 
Then $\tau$ is a coboundary  if and only if 
$\tau_a(g)=1$ for all $a\in A$ and $g\in G$ with $g\cdot a =a$. 
\end{Prop}
\begin{proof}
Assume that $\tau_a(g)=1$ for all $a\in A$ and $g\in G$ with $g\cdot a =a$. 
Fix representatives $\{a_1,\cdots,a_N\}$ of $A/G$. 
For $a=g\cdot a_l\in A$ we put $k_{a}:=\tau_{a_l}(g)\in K$. 
The cochain $(k_a)\in K^A$ is well-defined because 
for $g_1,g_2\in G$ with $g_1\cdot a=g_2\cdot a$, 
we have 
\begin{eqnarray*}
1=\tau_a(g_1-g_2)&=&\tau_{g_1\cdot a}(-g_2)\tau_a(g_1)\\
&=&\tau_{g_2\cdot a}(-g_2)\tau_a(g_1)\\
&=&\tau_a(g_2)^{-1}\tau_a(g_1).
\end{eqnarray*}
By definition we have 
$\tau_a(g)=k_{g\cdot a}k_a^{-1}$, and hence, 
$\tau$ is a coboundary.  
It is easy to check the converse. 
\end{proof}

\begin{proof}[Proof of Theorem~\ref{twistcohom}]
Define a homomorphism $\varphi:Z^1(\Gamma;\vec{j'})\to
\bigoplus_{[\vec j]\in \overline{QCG_k}(\Gamma;\vec{j'})}
H_1(\Gamma)_{\vec{j}}^*$ by 
$\varphi(\delta):\lambda\mapsto \delta_{\vec j}(\lambda)$ for 
$\delta\in Z^1(\Gamma;\vec{j'})$ and $\lambda\in H_1(\Gamma)_{\vec{j}}$. 
Note that the restriction 
$\delta_{\vec{j}}|_{H_1(\Gamma)_{\vec j}}$ depends only on 
$[{\vec j}]\in \overline{QCG_k}(\Gamma;\vec{j'})$ 
because of the relation 
\begin{eqnarray*}
\delta_{\lambda'\cdot\vec{j}}(\lambda)\delta_{\vec j}(\lambda')
&=&\delta_{\vec j}(\lambda+\lambda')\\
&=&\delta_{\lambda\cdot\vec{j}}(\lambda')
\delta_{\vec j}(\lambda)=
\delta_{\vec{j}}(\lambda')
\delta_{\vec j}(\lambda)
\end{eqnarray*}
for $\lambda\in H_1(\Gamma)_{\vec j}$ and $\lambda'\in H_1(\Gamma)$. 
It is clear that the kernel of $\varphi$ is 
equal to the group of coboundaries by Lemma~\ref{coboundary}. 
We show that $\varphi$ is surjective. 
Fix $\varepsilon=(\varepsilon_{[\vec j]}) \in \bigoplus_{[\vec j]\in \overline{QCG_k}(\Gamma;\vec{j'})}
H_1(\Gamma)_{\vec{j}}^*$ and 
take a lift $\tilde\varepsilon_{[\vec j]}:H_1(\Gamma)\to R^{\times}$
of $\varepsilon_{[\vec j]}:H_1(\Gamma)_{\vec j}\to R^{\times}$ 
for each $[\vec j]$. Define a map 
$\delta:H_1(\Gamma)\to A(\Gamma;\vec{j'})$ by 
$\delta_{\vec j}(\lambda):=\tilde\varepsilon_{[\vec j]}(\lambda)$. 
Then one can check that this $\delta$ is a twisted 
homomorphism and $\varphi(\delta)=\varepsilon$, 
and hence $\varphi$ is surjective. 
\end{proof}

\section{External edge condition}

In this section we recall the definition of 
the external edge condition (\cite{hjfj}). 
We first recall the definition of $\lambda$-external 
(resp. internal) edges. 

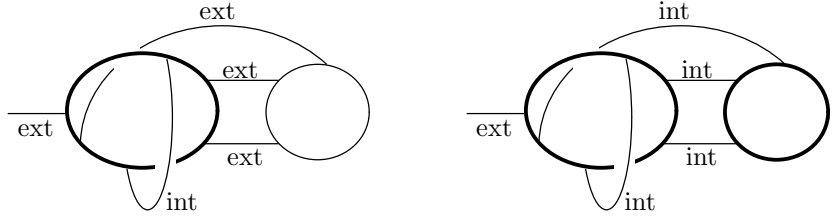
\begin{figure}[htb]
\input{externaledge.tex}
\caption{External edges and internal edges for $\lambda$ depicted by thick lines.}
\label{externaledge}
\end{figure}

\begin{Def}[$\lambda$-external/$\lambda$-internal edges]
Let $\Gamma$ be a unitrivalent graph. 
For a cycle $\lambda\in H_1(\Gamma)$, an edge $f_l\in \Gamma$ is said to be a 
$\lambda$-{\it external edge} if 
the cycle $\lambda$ on $\Gamma$ does not pass through $f_l$ 
and one of the vertex of $f_l$ lies on $\lambda$ and
the other is not. 
If $\lambda$ does not pass through $f_l$ and 
all vertices of $f_l$ lie on $\lambda$, then  
$f_l$ is said to be a {\it $\lambda$-internal edge}. 
See Figure~\ref{externaledge} for example. 
For $\lambda\in H_1(\Gamma)$ we denote the set of all 
$\lambda$-external edges by 
${\rm Ex}(\lambda)$.  
\end{Def}

\begin{Def}[External edge condition]\label{extedgecond}
Let $(\Gamma;\vec j')$ be a colored graph. 
Let $\delta=(\delta_{\vec j}):H_1(\Gamma)\to A(\Gamma;\vec j')$ 
be a twisted 1-cocycle. 
We say that $\delta$ satisfies the {\it external edge condition}
if the following condition is satisfied: 
$$
\delta_{\vec j}(\lambda)
=\left(-1\right)^{\sum_{f_l\in {\rm Ex}(\lambda)}j_l} \ 
{\rm for \ all} \ \lambda\in H_1(\Gamma) \ 
{\rm and} \ \vec j\in QCG_k(\Gamma;\vec j') \ 
{\rm with} \ \lambda\cdot \vec j=\vec j .  
$$
We call a twisted 1-cocycle satisfying the external edge condition 
an {\it external edge cocycle} and 
its cohomology class an {\it external edge class}. 
\end{Def}

\begin{Rem}
If $k$ is an odd number, then we have 
$H_1(\Gamma)_{\vec j}=\{0\}$ and 
$H^1(\Gamma;\vec j')=\{0\}$, and 
the trivial cocycle satisfies the external edge condition. 
Hereafter we will assume that $k$ is an even number. 
\end{Rem}

If we take $R$ to be the coefficient ring used in \cite{BHMV} 
and fix a ribbon structure of $\Gamma$, 
then there exists a canonical external edge cocycle as in 
\cite[Proposition~3.4]{hjfj}. 
Lemma~\ref{coboundary} implies its cohomology 
class does not depend on the ribbon structure. 
For general $R$ we have the following. 

\begin{Thm}\label{existextedge}
For each colored graph $(\Gamma;\vec j')$ 
there exist the unique external edge class. 
\end{Thm}

\begin{proof}
The uniqueness follows from Lemma~\ref{coboundary} 
and the definition of the external edge condition. 
By Theorem~\ref{twistcohom} it is enough to check that 
the map $\varepsilon=(\varepsilon_{[\vec j]}):
\bigoplus_{[\vec j]}H_1(\Gamma)_{\vec j}\to R^{\times}$ 
provided by the external edge condition, 
$$
\varepsilon_{[\vec j]}:\lambda\mapsto
\left(-1\right)^{\sum_{f_l\in {\rm Ex}(\lambda)}j_l}
\quad 
(\lambda\in H_1(\Gamma)_{\vec j}), 
$$is a group homomorphism. In other words 
we have to check the equality  
$$
\sum_{f_l\in {\rm Ex}(\lambda_1+\lambda_2)}j_l\equiv
\sum_{f_{l'}\in {\rm Ex}(\lambda_1)}j_{l'}+
\sum_{f_{l''}\in {\rm Ex}(\lambda_2)}j_{l''} \bmod 2 
$$for $\lambda_1,\lambda_2\in H_1(\Gamma)_{\vec j}$. 
Note that $\varepsilon_{[\vec j]}$ depends only on $[\vec j]$. 
In fact we have
$$
\sum_{f_l\in {\rm Ex}(\lambda)}\left((\lambda'\cdot\vec j)_l-
j_l\right)=
\sum_{f_{l}\in\lambda'\cap {\rm Ex}(\lambda)}\left(\frac{k}{2}-2j_{l}\right)
\equiv 0\bmod 2
$$ for all $\lambda'\in H_1(\Gamma)$ 
because 
$^{\#}(\lambda'\cap{\rm Ex}(\lambda))$ is even 
and the colorings 
on ${\rm Ex}(\lambda)$ are integers because of the $QCG_k$-condition. 
To show the required relation it is enough to 
consider two local 
situations in Figure~\ref{externalcocycle}.
\begin{figure}[htb]
\centering
\input{externalcocycle.tex}
\caption{Two local situations.}\label{externalcocycle}
\end{figure}
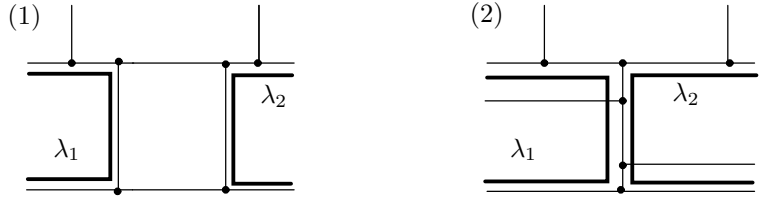
\\ By direct computations 
one can check that cancellations 
occur in these local situations. 
We only note that if the second 
situation occurs then we have $k\equiv 0 \bmod 4$ by 
the $QCG_k$-condition. 
\end{proof}

\section{Characterization of the external edge condition}

In this section we give a characterization of the 
external edge classes in a suitable categorical setting. 
We first introduce a representation associated with a twisted 1-cocycle, 
which gives a combinatorial realization of the Heisenberg representation 
studied in \cite{AndMas}, \cite{BHMV} and \cite{hjfj}. 
\begin{Def}\label{H1action}
For a colored graph $(\Gamma;\vec j')$ 
we denote by $R(\Gamma;\vec j')$ 
the free $R$-module generated by the finite set 
$QCG_k(\Gamma;\vec j')$ and by $|\vec j\rangle\in R(\Gamma;\vec j')$ 
the basis element corresponding to $\vec j\in QCG_k(\Gamma;\vec j')$. 
For a twisted 1-cocycle $\delta\in Z^1(\Gamma;\vec j')$ 
we can define a homomorphism 
$\rho(\delta):H_1(\Gamma)\to GL(R(\Gamma;\vec j'))$ 
by $\rho(\delta)(\lambda)|\vec j\rangle:=
\delta_{\vec j}(\lambda)|\lambda\cdot\vec j\rangle$.
\end{Def}

\begin{Prop}\label{coboundisom}
Let $\delta_1,\delta_2:H_1(\Gamma)\to A(\Gamma;\vec{j'})$ 
be two twisted 1-cocycles. 
If their cohomology classes coincide, 
$[\delta_1]=[\delta_2]\in H^1(\Gamma;\vec{j'})$,  then 
two representations $\rho(\delta_1)$ and $\rho(\delta_2)$ 
are isomorphic to each other. 
\end{Prop}

\begin{proof}
Fix $c=(c_{\vec j})\in C^0(\Gamma;\vec j')=A(\Gamma;\vec j')$ 
which cobounds $\delta_1\delta_2^{-1}$, i.e., 
$\delta_1(\lambda)\delta_2(\lambda)^{-1}=(\lambda\cdot c)c^{-1}$ 
for all $\lambda\in H_1(\Gamma)$. 
Then one can check that 
the map $\phi_c:R(\Gamma;\vec j')\to R(\Gamma;\vec j')$ defined by 
$\phi_c|\vec j\rangle=c_{\vec j}|\vec j\rangle$ gives 
an intertwiner between $\rho(\delta_1)$ and $\rho(\delta_2)$. 
\end{proof}

\begin{Rem}
By the construction we have $\phi_{c_1c_2}=\phi_{c_1}\phi_{c_2}$ 
for $c_1,c_2\in C^0(\Gamma;\vec j')$. 
This implies that the representation space 
$(R(\Gamma;\vec j'),\rho(\delta))$ is canonically defined for 
the cohomology class $[\delta]$. 
\end{Rem}

\begin{Rem}
The $R$-module $R(\Gamma;\vec j')$ is abstractly isomorphic to the 
TQFT-module in \cite{BHMV}, and 
the representation $\rho(\delta)$ gives a combinatorial realization of 
the Heisenberg representation on it. 
More precisely 
we can define a Heisenberg type group ${\mathcal E}(\Gamma)$ 
and an action of it on $R(\Gamma;\vec j')$ by using $\rho(\delta)$ . 
Note that if we fix a ribbon structure of $\Gamma$, 
then we have a colored surface $(C_{\Gamma};\vec j')$, and 
${\mathcal E}(\Gamma)$ is a central extension of $H_1(C_{\Gamma})$. 
See Section~5 of \cite{hjfj} for details. 
\cite[Proposition~3.4]{hjfj} and 
Proposition~\ref{coboundisom} above implies that 
if $\delta$ satisfies the external edge condition, 
then the induced Heisenberg representation defined by $\rho(\delta)$ 
is isomorphic to that studied in \cite{BHMV} and \cite{hjfj}. 
\end{Rem}

Consider a decomposition of the graph 
$\Gamma$ into $\Gamma_1\sqcup\Gamma_2$. 
In this article we mean a {\it decomposition of the graph} 
by cutting some internal edges in $\Gamma$. 
See Figure~\ref{decomposition} below.
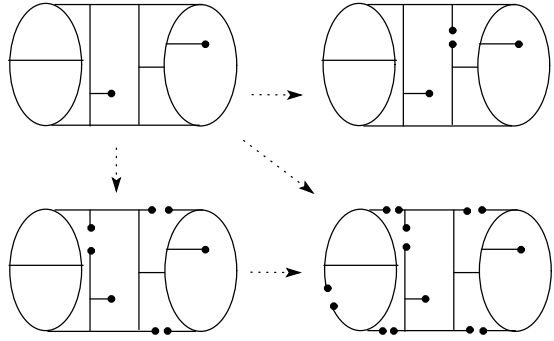
\begin{figure}[htb]
\begin{center}
\input{decomposition.tex}
\caption{Decompositions of the graph.}\label{decomposition}
\end{center}
\end{figure}
(Here $\Gamma_1$ or $\Gamma_2$ may be empty.)

Associated with the decomposition we have 
the decomposition 
$
QCG_k(\Gamma;\vec{j'})=\bigsqcup_{\vec{j''}}
QCG_k(\Gamma_1;\vec{j'_1})
\times QCG_k(\Gamma_2;\vec{j'_2})$, 
where $\vec{j}''$ denotes all colorings of the internal edges 
which are used to decompose $\Gamma$ and 
$\vec j'_1$ (resp. $\vec j'_2$) is the coloring 
which consists of values of $\vec j'\cup\vec{j}''$ on the 
univalent vertices on $\Gamma_1$ (resp.$\Gamma_2$). 
(Note that $\vec j'_1$ and $\vec j'_2$ 
have the common values on the components 
which come from $\vec{j}''$.) 
In this operation we say that 
$(\Gamma_1;\vec{j'_1})$ and 
$(\Gamma_2;\vec{j'_2})$ are a decomposition of 
$(\Gamma;\vec j')$, and we call 
$(\Gamma_2;\vec{j'_2})$ the complement of
$(\Gamma_1;\vec{j'_1})$ in $(\Gamma;\vec j')$. 
For each $\vec{j}\in QCG_k(\Gamma_2;\vec{j'_2})$ 
we have an inclusion 
$$\iota_{\vec j}:
QCG_k(\Gamma_1;\vec{j'_1})\times \{\vec{j}\}\hookrightarrow
QCG_k(\Gamma;\vec{j'})
$$and the induced projection 
$$
\iota_{\vec j}^*:A(\Gamma;\vec j')\to A(\Gamma_1;\vec{j_1'}).
$$
Combining with the inclusion 
$H_1(\Gamma_1)\hookrightarrow H_1(\Gamma)$ 
these data induce the cochain map 
and the map between cohomologies 
(denote them by the same letter)
$$
\iota_{\vec j}^*:C^*(\Gamma;\vec{j'})\to C^*(\Gamma_1;\vec{j_1'}), \ 
\iota_{\vec j}^*:H^*(\Gamma;\vec{j'})\to H^*(\Gamma_1;\vec{j_1'}).
$$

Now we define an equivalence relation $\sim_{\Gamma}$
of twisted 1-cocycles from the viewpoint of induced 
representations on the TQFT-module $R(\Gamma;\vec j')$. 

\begin{Def}
Let $(\Gamma;\vec j')$ be a colored graph.  
For two twisted 1-cocycles $\delta$ and $\delta'$ 
we denote by $\delta\sim_{\Gamma}\delta'$ 
if for any decomposition $(\Gamma_1;\vec{j'_1})$ and 
$(\Gamma_2;\vec{j'_2})$ 
of $(\Gamma;\vec j')$ 
the induced representations 
$\rho(\iota^*_{\vec j}(\delta))$ and $\rho(\iota^*_{\vec j}(\delta'))$ 
on $R(\Gamma_1;\vec j_1')$
are isomorphic representations 
for all $\vec j\in QCG_k(\Gamma_2;\vec j'_2)$. 
The relation $\sim_{\Gamma}$ is an equivalence relation. 
We denote by $[\delta]_{\Gamma}$ 
the equivalence class of a twisted 1-cocycle 
$\delta$ and  
by ${\mathcal H}(\Gamma;\vec j')$ the set of equivalence classes 
of twisted 1-cocycles under the equivalence 
relation $\sim_{\Gamma}$. 

Note that ${\mathcal H}(\Gamma;\vec j')$ has 
a structure of an abelian group 
induced from that of $Z^1(\Gamma;\vec j')$,  
$[\delta]_{\Gamma}\cdot[\delta']_{\Gamma}:=
[\delta\delta']_{\Gamma}$. 
\end{Def}

\begin{Thm}\label{isom&char}
The map 
$F:{\mathcal H}(\Gamma;\vec j')\to H^1(\Gamma;\vec j')$ 
defined by 
$F([\delta]_{\Gamma})=[\delta]$ is a group isomorphism.
\end{Thm}

To show this theorem we show the well-definedness and 
the injectivity of the map $F$. 
(The surjectivity is clear.) 
We first consider the following special case. 

\begin{Prop}\label{H_1=Z_2}
Let $(\Gamma(n),\vec j')$ be a colored graph, 
where $\Gamma(n)$ is a unitrivalent graph 
having $n$ univalent vertices and satisfying 
$H_1(\Gamma(n))\cong\Z/2$ (Figure~\ref{Gamman}). 
Let $\delta_1$ and $\delta_2$ be twisted 1-cocycles. 
Then two representations $\rho(\delta_1)$ and 
$\rho(\delta_2)$ are isomorphic if and only if 
their cohomology classes coincide, 
$[\delta_1]=[\delta_2]\in H^1(\Gamma(n);\vec j')$. 
\end{Prop}

\begin{figure}[htb]
\centering
\input{gamman.tex}
\caption{$\Gamma(n)$.}\label{gamman}
\end{figure}

\begin{proof}
Let $\lambda$ be the generator of $H_1(\Gamma(n))$. 
By Lemma~\ref{coboundary} it is enough to show that 
if $\rho(\delta_1)$ and 
$\rho(\delta_2)$ are isomorphic, then 
$(\delta_1)_{\vec j}(\lambda)=(\delta_2)_{\vec j}(\lambda)$ 
for all pairs $(\lambda,\vec j)$ with $\lambda\cdot\vec j=\vec j$. 
Note that if there is such $\vec j$, 
then it is equal to 
$\vec j_0:=(\vec j',\frac{k}{4},\cdots,\frac{k}{4})$. 
By definition of the representation $\rho(\delta_i)$ we have 
${\rm Tr}(\rho(\delta_i)(\lambda))=(\delta_i)_{\vec j_0}(\lambda)$ for $i=1,2$. 
Since $\rho(\delta_1)$ and 
$\rho(\delta_2)$ are isomorphic, 
we have 
$(\delta_1)_{\vec j_0}(\lambda)=
(\delta_2)_{\vec j_0}(\lambda)$. 
\end{proof}

\begin{Lem}\label{well-def}
The map $F$ is well-defined. 
\end{Lem}
\begin{proof}
By Lemma~\ref{coboundary} it is enough to show that 
if $\delta_1\sim_{\Gamma}\delta_2$, 
then 
$(\delta_1)_{\vec j}(\lambda)=(\delta_2)_{\vec j}(\lambda)$ 
for all pairs $(\lambda,\vec j)$ with $\lambda\cdot\vec j=\vec j$. 
Fix such a pair $(\lambda,\vec j)$.
Take a decomposition 
$(\Gamma(\lambda);\vec j'_1)$ and $(\Gamma'(\lambda);\vec j'_2)$ 
of $(\Gamma;\vec j')$ as follows. 
We first 
cut all $\lambda$-external and $\lambda$-internal edges of $\Gamma$. 
Then we define 
$\Gamma(\lambda)$ to be the component of this graph 
which contains the cycle $\lambda$ and 
$\Gamma'(\lambda)$ to be the complementary graph. 
(See Figure~\ref{Gammalambda}.)
The colorings $\vec j'_1$ and $\vec j'_2$ 
are the induced colorings from $\vec j'$ and $\vec j$. 
\begin{figure}[htb]
\begin{center}
\input{gammalambda.tex}
\caption{$\Gamma(\lambda)$ and $\Gamma'(\lambda)$.}\label{Gammalambda}
\end{center}
\end{figure}
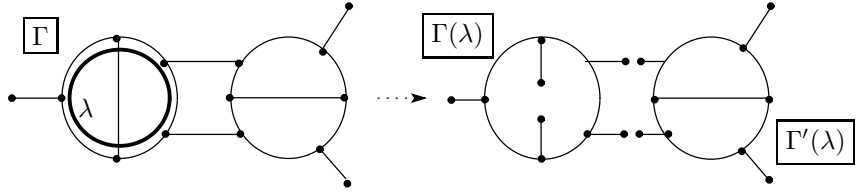
Note that $\Gamma(\lambda)$ is equal to 
$\Gamma(n')$ for some non-negative integer 
$n'$ (or as a disjoint union of such graphs). 
By the assumption induced representations 
$\rho(\iota_{\vec j_2}^*(\delta_1))$ and 
$\rho(\iota_{\vec j_2}^*(\delta_2))$ on 
$R(\Gamma(\lambda);\vec j'_1)$
are isomorphic representations 
for all $\vec j_2\in QCG_k(\Gamma'(\lambda);\vec j_2')$. 
In particular by Proposition~\ref{H_1=Z_2} 
we have 
$
(\delta_1)_{\vec j}(\lambda)=
\left(\iota_{\vec j_2}^*(\delta_1)\right)_{\vec j_1}(\lambda)=
\left(\iota_{\vec j_2}^*(\delta_2)\right)_{\vec j_1}(\lambda)
=(\delta_2)_{\vec j}(\lambda),
$ 
where $\vec j_1\in QCG_k(\Gamma(\lambda);\vec j_1')$ 
satisfies $\iota_{\vec j_2}(\vec j_1)=\vec j$ 
and $\lambda\cdot\vec j_1=\vec j_1$.
\end{proof}

\begin{Lem}
The map $F$ is injective. 
\end{Lem}
\begin{proof}
We show that 
if $F([\delta]_{\Gamma})=[\delta]=0$, 
then $\delta$ is equivalent to $0$. 
Take a decomposition $(\Gamma_1;\vec{j'_1})$ and 
$(\Gamma_2;\vec{j'_2})$
of $(\Gamma;\vec j')$. 
We have a homomorphism 
$\iota_{\vec j}^*:H^1(\Gamma;\vec{j'})\to H^1(\Gamma_1;\vec{j_1'})$ 
for each $\vec j\in QCG_k(\Gamma_2;\vec j'_2)$. 
By the assumption we have 
$\iota^*_{\vec j}[\delta]=0$, and hence the induced representation 
$\rho(\iota^*_{\vec j}(\delta))$ on $R(\Gamma_1;\vec{j_1'})$ 
is isomorphic to $\rho(0)$ by Proposition~\ref{coboundisom}. 
In particular we have $\delta\sim_{\Gamma}0$.
\end{proof}

We define two categories to state the naturality 
and a characterization of the external edge class.
\begin{Def}
For a positive integer $k$ 
we define a category ${\mathcal G}_k$ as follows. 
Objects of ${\mathcal G}_k$ consist of all 
colored graphs $\{(\Gamma;\vec j')\}$. 
The set of morphisms between two colored graphs $(\Gamma_1;\vec j'_1)$ 
and $(\Gamma_2;\vec j'_2)$ is empty if 
$(\Gamma_2; \vec j'_2)$ can not be obtained as any 
decompositions of $(\Gamma_1;\vec j'_1)$. 
If $(\Gamma_2;\vec j'_2)$ is a  decomposition of $(\Gamma_1;\vec j'_1)$ 
with the complement $(\Gamma_3;\vec j'_3)$, then 
we define the set of morphisms by 
${\rm Mor}((\Gamma_1;\vec j'_1),(\Gamma_2;\vec j'_2))
:=QCG_k(\Gamma_3;\vec j'_3)$. 
The composition of morphisms for a triple 
$(\Gamma_1;\vec j'_1)$, $(\Gamma_2;\vec j'_2)$ and 
$(\Gamma_3;\vec j'_3)$
is defined by the natural map 
$QCG_k(\Gamma_4;\vec j'_4)\times QCG_k(\Gamma_5;\vec j'_5)\to 
QCG_k(\Gamma_6;\vec j'_6)$, where 
$(\Gamma_4;\vec j'_4)$ is the complement of 
$(\Gamma_2;\vec j'_2)$ in $(\Gamma_1;\vec j'_1)$, 
$(\Gamma_5;\vec j'_5)$ is the complement of 
$(\Gamma_3;\vec j'_3)$ in $(\Gamma_2;\vec j'_2)$ and 
$(\Gamma_6;\vec j'_6)$ is the complement of 
$(\Gamma_3;\vec j'_3)$ in $(\Gamma_1;\vec j'_1)$. 

Let ${\mathcal Vec}^{0}$ be the category 
whose objects consist of pairs of a vector space 
$V$ and a vector $v\in V$ 
and morphisms consist of linear maps preserving the prescribed vectors. 
\end{Def}

As we showed in Theorem~\ref{existextedge} 
there exists the unique external edge class 
for each $(\Gamma;\vec j')$, and  
we denote it by $\delta_{\rm ext}[\Gamma;\vec j']\in H^1(\Gamma;\vec j')$. 
By definition we have the following. 

\begin{Prop}\label{universal}
We can define a functor 
${\mathcal  F}_{\rm ext}:{\mathcal G}_k\to {\mathcal Vec}^0$ 
by putting 
${\mathcal F}_{\rm ext}(\Gamma;\vec j')=
(H^1(\Gamma;\vec j'), \delta_{\rm ext}[\Gamma;\vec j'])$ and 
${\mathcal F}_{\rm ext}(\vec j)=\iota^*_{\vec j}$ for 
a decomposition $(\Gamma_1;\vec j'_1)$ and 
$(\Gamma_2;\vec j'_2)$ of $(\Gamma;\vec j')$ and 
$\vec j\in QCG_k(\Gamma_2;\vec j'_2)$. 
\end{Prop}

\begin{Rem}
By the similar argument in the proof of Lemma~\ref{well-def} 
we can check that a functor 
${\mathcal F}:{\mathcal G}_k\to {\mathcal Vec}^0$ 
of the form ${\mathcal F}(\Gamma;\vec j')=(H^1(\Gamma;\vec j'), 
\delta[\Gamma;\vec j'])$ and ${\mathcal F}(\vec j)=\iota^*_{\vec j}$ 
is characterized by values for all $(\Gamma(n),\vec j')$. 
\end{Rem}

As a corollary of Theorem~\ref{isom&char} 
we have the following characterization. 

\begin{Thm}\label{charofext}
For $l=1,2$ let 
${\mathcal F}_l:{\mathcal G}_k\to {\mathcal Vec}^0$ be two functors 
of the form ${\mathcal F}_l(\Gamma;\vec j')=(H^1(\Gamma;\vec j'), 
\delta_l[\Gamma;\vec j'])$ and ${\mathcal F}_l(\vec j)=\iota^*_{\vec j}$.  
If the induced representations $\rho(\delta_1[\Gamma;\vec j'])$ and 
$\rho(\delta_2[\Gamma;\vec j'])$ of $H_1(\Gamma)$ 
on $R(\Gamma;\vec j')$ are isomorphic to 
each other for all $(\Gamma;\vec j')$, 
then we have ${\mathcal F}_1={\mathcal F}_2$. 
In particular the functor ${\mathcal F}_{\rm ext}$ is the 
unique functor which realizes the representations studied in 
\cite{AndMas}, \cite{BHMV} and \cite{hjfj}. 
\end{Thm}

\end{document}

%% file: h1action.tex
\unitlength 0.1in
\begin{picture}( 19.9300,  7.0300)(  8.8000,-14.6000)
%
\special{pn 13}%
\special{ar 1484 1134 322 322  0.0000000 6.2831853}%
\put(26.9000,-12.1000){\makebox(0,0)[lb]{10}}%
\put(22.9100,-12.1100){\makebox(0,0)[lb]{9}}%
\put(19.6300,-14.1200){\makebox(0,0)[lb]{8}}%
\put(19.6300,-9.2700){\makebox(0,0)[lb]{7}}%
\put(17.1300,-12.0100){\makebox(0,0)[lb]{6}}%
\put(13.0600,-13.9200){\makebox(0,0)[lb]{5}}%
\put(12.0300,-11.8200){\makebox(0,0)[lb]{4}}%
\put(13.3100,-9.5100){\makebox(0,0)[lb]{3}}%
\put(9.7800,-14.9000){\makebox(0,0)[lb]{2}}%
\put(9.5300,-9.9600){\makebox(0,0)[lb]{1}}%
%
\special{pn 20}%
\special{sh 1}%
\special{ar 886 1418 10 10 0  6.28318530717959E+0000}%
%
\special{pn 20}%
\special{sh 1}%
\special{ar 1204 1286 10 10 0  6.28318530717959E+0000}%
%
\special{pn 20}%
\special{sh 1}%
\special{ar 1214 966 10 10 0  6.28318530717959E+0000}%
%
\special{pn 20}%
\special{sh 1}%
\special{ar 1478 814 10 10 0  6.28318530717959E+0000}%
%
\special{pn 20}%
\special{sh 1}%
\special{ar 1458 1456 10 10 0  6.28318530717959E+0000}%
%
\special{pn 20}%
\special{sh 1}%
\special{ar 2542 1456 10 10 0  6.28318530717959E+0000}%
%
\special{pn 20}%
\special{sh 1}%
\special{ar 2556 820 10 10 0  6.28318530717959E+0000}%
%
\special{pn 20}%
\special{sh 1}%
\special{ar 890 830 10 10 0  6.28318530717959E+0000}%
%
\special{pn 8}%
\special{pa 1194 1290}%
\special{pa 880 1412}%
\special{fp}%
%
\special{pn 8}%
\special{pa 1208 972}%
\special{pa 890 824}%
\special{fp}%
%
\special{pn 8}%
\special{ar 2552 1138 322 322  0.0000000 6.2831853}%
%
\special{pn 8}%
\special{pa 1478 1456}%
\special{pa 2556 1456}%
\special{fp}%
%
\special{pn 8}%
\special{pa 1474 810}%
\special{pa 2552 810}%
\special{fp}%
\end{picture}%

%% file: externaledge.tex
\unitlength 0.1in
\begin{picture}( 42.9000, 11.6600)( 10.3800,-16.0200)
%
\special{pn 8}%
\special{ar 2160 1270 744 632  3.1681016 5.5323330}%
%
\special{pn 8}%
\special{pa 1348 1098}%
\special{pa 1038 1098}%
\special{fp}%
\put(10.8800,-12.1400){\makebox(0,0)[lb]{ext}}%
%
\special{pn 8}%
\special{pa 2068 924}%
\special{pa 2452 924}%
\special{fp}%
%
\special{pn 8}%
\special{pa 2068 1256}%
\special{pa 2452 1256}%
\special{fp}%
%
\special{pn 8}%
\special{ar 2658 1090 270 250  0.1354594 6.2831853}%
\special{ar 2658 1090 270 250  0.0000000 0.1314492}%
\put(18.6400,-16.0600){\makebox(0,0)[lb]{int}}%
%
\special{pn 8}%
\special{sh 0}%
\special{pa 1816 1142}%
\special{pa 1924 1142}%
\special{pa 1924 1250}%
\special{pa 1816 1250}%
\special{pa 1816 1142}%
\special{ip}%
\put(21.8400,-13.6600){\makebox(0,0)[lb]{ext}}%
\put(21.6000,-9.1000){\makebox(0,0)[lb]{ext}}%
%
\special{pn 8}%
\special{sh 0}%
\special{pa 1592 734}%
\special{pa 1728 734}%
\special{pa 1728 910}%
\special{pa 1592 910}%
\special{pa 1592 734}%
\special{ip}%
%
\special{pn 20}%
\special{ar 1740 1082 394 300  0.0079680 6.2831853}%
%
\special{pn 8}%
\special{sh 0}%
\special{pa 1808 1326}%
\special{pa 1916 1326}%
\special{pa 1916 1434}%
\special{pa 1808 1434}%
\special{pa 1808 1326}%
\special{ip}%
%
\special{pn 8}%
\special{ar 1768 1118 130 484  5.6039161 6.2831853}%
\special{ar 1768 1118 130 484  0.0000000 2.5467584}%
\put(20.4000,-6.0600){\makebox(0,0)[lb]{ext}}%
%
\special{pn 8}%
\special{ar 4560 1270 744 632  3.1681016 5.5323330}%
%
\special{pn 8}%
\special{pa 3750 1098}%
\special{pa 3440 1098}%
\special{fp}%
\put(34.8800,-12.1400){\makebox(0,0)[lb]{ext}}%
%
\special{pn 8}%
\special{pa 4470 924}%
\special{pa 4854 924}%
\special{fp}%
%
\special{pn 8}%
\special{pa 4470 1256}%
\special{pa 4854 1256}%
\special{fp}%
%
\special{pn 20}%
\special{ar 5060 1090 270 250  0.1354594 6.2831853}%
\special{ar 5060 1090 270 250  0.0000000 0.1314492}%
\put(42.6500,-16.0600){\makebox(0,0)[lb]{int}}%
%
\special{pn 8}%
\special{sh 0}%
\special{pa 4218 1142}%
\special{pa 4324 1142}%
\special{pa 4324 1250}%
\special{pa 4218 1250}%
\special{pa 4218 1142}%
\special{ip}%
\put(45.8500,-13.6600){\makebox(0,0)[lb]{int}}%
\put(45.6100,-9.1000){\makebox(0,0)[lb]{int}}%
%
\special{pn 8}%
\special{sh 0}%
\special{pa 3994 734}%
\special{pa 4130 734}%
\special{pa 4130 910}%
\special{pa 3994 910}%
\special{pa 3994 734}%
\special{ip}%
%
\special{pn 20}%
\special{ar 4142 1082 394 300  0.0079759 6.2831853}%
%
\special{pn 8}%
\special{sh 0}%
\special{pa 4210 1326}%
\special{pa 4316 1326}%
\special{pa 4316 1434}%
\special{pa 4210 1434}%
\special{pa 4210 1326}%
\special{ip}%
%
\special{pn 8}%
\special{ar 4170 1118 132 484  5.6031970 6.2831853}%
\special{ar 4170 1118 132 484  0.0000000 2.5442127}%
\put(44.4100,-6.0600){\makebox(0,0)[lb]{int}}%
\end{picture}%

%% file: externalcocycle.tex
\unitlength 0.1in
\begin{picture}( 39.2000, 10.1600)(  6.9000,-19.0600)
%
\special{pn 8}%
\special{pa 796 1232}%
\special{pa 1350 1232}%
\special{fp}%
%
\special{pn 8}%
\special{pa 1270 1232}%
\special{pa 1270 1890}%
\special{fp}%
%
\special{pn 8}%
\special{pa 790 1896}%
\special{pa 1344 1896}%
\special{fp}%
%
\special{pn 8}%
\special{pa 1344 1896}%
\special{pa 2190 1896}%
\special{fp}%
%
\special{pn 8}%
\special{pa 1344 1232}%
\special{pa 2190 1232}%
\special{fp}%
%
\special{pn 8}%
\special{pa 1832 1238}%
\special{pa 1832 1896}%
\special{fp}%
%
\special{pn 8}%
\special{pa 1026 1232}%
\special{pa 1026 930}%
\special{fp}%
%
\special{pn 8}%
\special{pa 2006 1232}%
\special{pa 2006 930}%
\special{fp}%
%
\special{pn 20}%
\special{pa 802 1288}%
\special{pa 1228 1288}%
\special{pa 1228 1840}%
\special{pa 796 1840}%
\special{pa 796 1840}%
\special{fp}%
%
\special{pn 20}%
\special{sh 1}%
\special{ar 1026 1232 10 10 0  6.28318530717959E+0000}%
%
\special{pn 20}%
\special{sh 1}%
\special{ar 2000 1232 10 10 0  6.28318530717959E+0000}%
%
\special{pn 20}%
\special{sh 1}%
\special{ar 1832 1238 10 10 0  6.28318530717959E+0000}%
%
\special{pn 20}%
\special{sh 1}%
\special{ar 1832 1896 10 10 0  6.28318530717959E+0000}%
%
\special{pn 20}%
\special{sh 1}%
\special{ar 1266 1904 10 10 0  6.28318530717959E+0000}%
%
\special{pn 20}%
\special{sh 1}%
\special{ar 1272 1224 10 10 0  6.28318530717959E+0000}%
%
\special{pn 20}%
\special{pa 2180 1296}%
\special{pa 1872 1296}%
\special{pa 1872 1864}%
\special{pa 2174 1864}%
\special{pa 2174 1864}%
\special{fp}%
\put(6.9000,-10.6000){\makebox(0,0)[lb]{(1)}}%
\put(31.1000,-10.5000){\makebox(0,0)[lb]{(2)}}%
\put(9.3000,-17.4300){\makebox(0,0)[lb]{$\lambda_1$}}%
\put(20.1400,-14.6100){\makebox(0,0)[lb]{$\lambda_2$}}%
%
\special{pn 8}%
\special{pa 3200 1232}%
\special{pa 4600 1232}%
\special{fp}%
%
\special{pn 8}%
\special{pa 3200 1904}%
\special{pa 4600 1904}%
\special{fp}%
%
\special{pn 8}%
\special{pa 3910 1232}%
\special{pa 3910 1904}%
\special{fp}%
%
\special{pn 8}%
\special{pa 2006 1232}%
\special{pa 2006 930}%
\special{fp}%
%
\special{pn 8}%
\special{pa 3500 1232}%
\special{pa 3500 930}%
\special{fp}%
%
\special{pn 8}%
\special{pa 4460 1232}%
\special{pa 4460 930}%
\special{fp}%
%
\special{pn 8}%
\special{pa 3910 1430}%
\special{pa 3190 1430}%
\special{fp}%
\special{pa 3920 1762}%
\special{pa 4600 1762}%
\special{fp}%
%
\special{pn 20}%
\special{sh 1}%
\special{ar 4470 1232 10 10 0  6.28318530717959E+0000}%
%
\special{pn 20}%
\special{sh 1}%
\special{ar 3900 1896 10 10 0  6.28318530717959E+0000}%
%
\special{pn 20}%
\special{sh 1}%
\special{ar 3910 1768 10 10 0  6.28318530717959E+0000}%
%
\special{pn 20}%
\special{sh 1}%
\special{ar 3910 1430 10 10 0  6.28318530717959E+0000}%
%
\special{pn 20}%
\special{sh 1}%
\special{ar 3910 1232 10 10 0  6.28318530717959E+0000}%
%
\special{pn 20}%
\special{sh 1}%
\special{ar 3500 1232 10 10 0  6.28318530717959E+0000}%
%
\special{pn 20}%
\special{pa 4610 1296}%
\special{pa 3960 1296}%
\special{pa 3960 1858}%
\special{pa 4590 1858}%
\special{pa 4590 1858}%
\special{fp}%
%
\special{pn 20}%
\special{pa 3200 1308}%
\special{pa 3830 1308}%
\special{pa 3830 1852}%
\special{pa 3190 1852}%
\special{pa 3190 1852}%
\special{fp}%
\put(33.2000,-17.4300){\makebox(0,0)[lb]{$\lambda_1$}}%
\put(41.7000,-14.5500){\makebox(0,0)[lb]{$\lambda_2$}}%
\end{picture}%

%% file: decomposition.tex
\unitlength 0.1in
\begin{picture}( 28.4000, 17.8500)( 28.7800,-31.9900)
%
\special{pn 8}%
\special{pa 5200 1420}%
\special{pa 5200 2048}%
\special{fp}%
%
\special{pn 8}%
\special{sh 0}%
\special{pa 5178 1554}%
\special{pa 5228 1554}%
\special{pa 5228 1640}%
\special{pa 5178 1640}%
\special{pa 5178 1554}%
\special{ip}%
%
\special{pn 8}%
\special{pa 4952 2962}%
\special{pa 5070 2962}%
\special{fp}%
\special{pa 5208 2824}%
\special{pa 5340 2824}%
\special{fp}%
%
\special{pn 8}%
\special{pa 3120 2496}%
\special{pa 3894 2496}%
\special{fp}%
%
\special{pn 8}%
\special{sh 0}%
\special{pa 3636 2464}%
\special{pa 3708 2464}%
\special{pa 3708 2546}%
\special{pa 3636 2546}%
\special{pa 3636 2464}%
\special{ip}%
%
\special{pn 8}%
\special{pa 4758 3128}%
\special{pa 5530 3128}%
\special{fp}%
%
\special{pn 8}%
\special{sh 0}%
\special{pa 5304 3114}%
\special{pa 5356 3114}%
\special{pa 5356 3200}%
\special{pa 5304 3200}%
\special{pa 5304 3114}%
\special{ip}%
%
\special{pn 8}%
\special{sh 0}%
\special{pa 4844 3110}%
\special{pa 4896 3110}%
\special{pa 4896 3194}%
\special{pa 4844 3194}%
\special{pa 4844 3110}%
\special{ip}%
%
\special{pn 8}%
\special{pa 4768 2496}%
\special{pa 5540 2496}%
\special{fp}%
%
\special{pn 8}%
\special{sh 0}%
\special{pa 5284 2464}%
\special{pa 5356 2464}%
\special{pa 5356 2546}%
\special{pa 5284 2546}%
\special{pa 5284 2464}%
\special{ip}%
%
\special{pn 8}%
\special{ar 5530 2814 188 318  0.0000000 6.2831853}%
%
\special{pn 20}%
\special{sh 1}%
\special{ar 5278 2502 10 10 0  6.28318530717959E+0000}%
%
\special{pn 8}%
\special{sh 0}%
\special{pa 4854 2450}%
\special{pa 4906 2450}%
\special{pa 4906 2536}%
\special{pa 4854 2536}%
\special{pa 4854 2450}%
\special{ip}%
%
\special{pn 8}%
\special{ar 4720 2810 188 318  0.0000000 6.2831853}%
%
\special{pn 8}%
\special{sh 0}%
\special{pa 4532 2908}%
\special{pa 4584 2908}%
\special{pa 4584 2994}%
\special{pa 4532 2994}%
\special{pa 4532 2908}%
\special{ip}%
%
\special{pn 8}%
\special{pa 3304 2962}%
\special{pa 3422 2962}%
\special{fp}%
\special{pa 3560 2824}%
\special{pa 3694 2824}%
\special{fp}%
%
\special{pn 8}%
\special{pa 3560 2496}%
\special{pa 3560 3122}%
\special{fp}%
%
\special{pn 8}%
\special{pa 3304 2496}%
\special{pa 3304 3132}%
\special{fp}%
%
\special{pn 8}%
\special{pa 3084 3132}%
\special{pa 3856 3132}%
\special{fp}%
%
\special{pn 8}%
\special{sh 0}%
\special{pa 3656 3114}%
\special{pa 3708 3114}%
\special{pa 3708 3200}%
\special{pa 3656 3200}%
\special{pa 3656 3114}%
\special{ip}%
%
\special{pn 8}%
\special{ar 3072 1734 188 320  0.0000000 6.2831853}%
%
\special{pn 8}%
\special{pa 4108 2132}%
\special{pa 4470 2406}%
\special{dt 0.045}%
\special{sh 1}%
\special{pa 4470 2406}%
\special{pa 4430 2350}%
\special{pa 4428 2374}%
\special{pa 4406 2382}%
\special{pa 4470 2406}%
\special{fp}%
%
\special{pn 8}%
\special{pa 4146 2822}%
\special{pa 4402 2824}%
\special{dt 0.045}%
\special{sh 1}%
\special{pa 4402 2824}%
\special{pa 4336 2804}%
\special{pa 4348 2824}%
\special{pa 4334 2844}%
\special{pa 4402 2824}%
\special{fp}%
%
\special{pn 8}%
\special{pa 4152 1892}%
\special{pa 4408 1894}%
\special{dt 0.045}%
\special{sh 1}%
\special{pa 4408 1894}%
\special{pa 4342 1874}%
\special{pa 4356 1894}%
\special{pa 4342 1914}%
\special{pa 4408 1894}%
\special{fp}%
%
\special{pn 8}%
\special{pa 3442 2172}%
\special{pa 3442 2396}%
\special{dt 0.045}%
\special{sh 1}%
\special{pa 3442 2396}%
\special{pa 3462 2330}%
\special{pa 3442 2344}%
\special{pa 3422 2330}%
\special{pa 3442 2396}%
\special{fp}%
%
\special{pn 20}%
\special{sh 1}%
\special{ar 5202 1556 10 10 0  6.28318530717959E+0000}%
%
\special{pn 20}%
\special{sh 1}%
\special{ar 5202 1628 10 10 0  6.28318530717959E+0000}%
%
\special{pn 20}%
\special{sh 1}%
\special{ar 5548 1628 10 10 0  6.28318530717959E+0000}%
%
\special{pn 20}%
\special{sh 1}%
\special{ar 3416 2962 10 10 0  6.28318530717959E+0000}%
%
\special{pn 20}%
\special{sh 1}%
\special{ar 3646 3130 10 10 0  6.28318530717959E+0000}%
%
\special{pn 20}%
\special{sh 1}%
\special{ar 3908 2704 10 10 0  6.28318530717959E+0000}%
%
\special{pn 20}%
\special{sh 1}%
\special{ar 3628 2496 10 10 0  6.28318530717959E+0000}%
%
\special{pn 20}%
\special{sh 1}%
\special{ar 3716 2496 10 10 0  6.28318530717959E+0000}%
%
\special{pn 20}%
\special{sh 1}%
\special{ar 4856 2496 10 10 0  6.28318530717959E+0000}%
%
\special{pn 20}%
\special{sh 1}%
\special{ar 4920 2496 10 10 0  6.28318530717959E+0000}%
%
\special{pn 20}%
\special{sh 1}%
\special{ar 5356 2496 10 10 0  6.28318530717959E+0000}%
%
\special{pn 20}%
\special{sh 1}%
\special{ar 4548 2906 10 10 0  6.28318530717959E+0000}%
%
\special{pn 20}%
\special{sh 1}%
\special{ar 4580 3000 10 10 0  6.28318530717959E+0000}%
%
\special{pn 20}%
\special{sh 1}%
\special{ar 4836 3130 10 10 0  6.28318530717959E+0000}%
%
\special{pn 20}%
\special{sh 1}%
\special{ar 4894 3130 10 10 0  6.28318530717959E+0000}%
%
\special{pn 20}%
\special{sh 1}%
\special{ar 5060 2962 10 10 0  6.28318530717959E+0000}%
%
\special{pn 20}%
\special{sh 1}%
\special{ar 5292 3124 10 10 0  6.28318530717959E+0000}%
%
\special{pn 20}%
\special{sh 1}%
\special{ar 5362 3130 10 10 0  6.28318530717959E+0000}%
%
\special{pn 20}%
\special{sh 1}%
\special{ar 5560 2704 10 10 0  6.28318530717959E+0000}%
%
\special{pn 20}%
\special{sh 1}%
\special{ar 5080 1886 10 10 0  6.28318530717959E+0000}%
%
\special{pn 20}%
\special{sh 1}%
\special{ar 3908 1628 10 10 0  6.28318530717959E+0000}%
%
\special{pn 20}%
\special{sh 1}%
\special{ar 3416 1886 10 10 0  6.28318530717959E+0000}%
%
\special{pn 8}%
\special{pa 5352 2702}%
\special{pa 5562 2702}%
\special{fp}%
%
\special{pn 8}%
\special{pa 4528 2786}%
\special{pa 4916 2786}%
\special{fp}%
%
\special{pn 8}%
\special{pa 4952 2496}%
\special{pa 4952 3132}%
\special{fp}%
%
\special{pn 8}%
\special{pa 5208 2496}%
\special{pa 5208 3122}%
\special{fp}%
%
\special{pn 8}%
\special{pa 5342 1626}%
\special{pa 5552 1626}%
\special{fp}%
%
\special{pn 8}%
\special{pa 4518 1712}%
\special{pa 4908 1712}%
\special{fp}%
%
\special{pn 8}%
\special{pa 4944 1420}%
\special{pa 4944 2056}%
\special{fp}%
%
\special{pn 8}%
\special{pa 4736 2056}%
\special{pa 5510 2056}%
\special{fp}%
%
\special{pn 8}%
\special{pa 4758 1420}%
\special{pa 5532 1420}%
\special{fp}%
%
\special{pn 8}%
\special{ar 5522 1738 188 318  0.0000000 6.2831853}%
%
\special{pn 8}%
\special{pa 4944 1886}%
\special{pa 5060 1886}%
\special{fp}%
\special{pa 5200 1748}%
\special{pa 5332 1748}%
\special{fp}%
%
\special{pn 8}%
\special{ar 4712 1734 188 320  0.0000000 6.2831853}%
%
\special{pn 8}%
\special{pa 3704 2702}%
\special{pa 3914 2702}%
\special{fp}%
%
\special{pn 8}%
\special{pa 2880 2786}%
\special{pa 3268 2786}%
\special{fp}%
%
\special{pn 8}%
\special{ar 3072 2810 188 318  0.0000000 6.2831853}%
%
\special{pn 8}%
\special{pa 3702 1626}%
\special{pa 3912 1626}%
\special{fp}%
%
\special{pn 8}%
\special{pa 2878 1712}%
\special{pa 3268 1712}%
\special{fp}%
%
\special{pn 8}%
\special{pa 3304 1420}%
\special{pa 3304 2056}%
\special{fp}%
%
\special{pn 8}%
\special{pa 3560 1420}%
\special{pa 3560 2048}%
\special{fp}%
%
\special{pn 8}%
\special{pa 3104 2052}%
\special{pa 3876 2052}%
\special{fp}%
%
\special{pn 8}%
\special{pa 3120 1420}%
\special{pa 3892 1420}%
\special{fp}%
%
\special{pn 8}%
\special{ar 3882 1738 190 318  0.0000000 6.2831853}%
%
\special{pn 8}%
\special{pa 3304 1886}%
\special{pa 3420 1886}%
\special{fp}%
\special{pa 3560 1748}%
\special{pa 3692 1748}%
\special{fp}%
%
\special{pn 20}%
\special{sh 1}%
\special{ar 3710 3130 10 10 0  6.28318530717959E+0000}%
%
\special{pn 8}%
\special{ar 3884 2814 188 318  0.0000000 6.2831853}%
%
\special{pn 8}%
\special{sh 0}%
\special{pa 3280 2580}%
\special{pa 3340 2580}%
\special{pa 3340 2700}%
\special{pa 3280 2700}%
\special{pa 3280 2580}%
\special{ip}%
%
\special{pn 20}%
\special{sh 1}%
\special{ar 3310 2590 10 10 0  6.28318530717959E+0000}%
%
\special{pn 20}%
\special{sh 1}%
\special{ar 3310 2710 10 10 0  6.28318530717959E+0000}%
%
\special{pn 8}%
\special{sh 0}%
\special{pa 4920 2580}%
\special{pa 4980 2580}%
\special{pa 4980 2700}%
\special{pa 4920 2700}%
\special{pa 4920 2580}%
\special{ip}%
%
\special{pn 20}%
\special{sh 1}%
\special{ar 4960 2590 10 10 0  6.28318530717959E+0000}%
%
\special{pn 20}%
\special{sh 1}%
\special{ar 4960 2690 10 10 0  6.28318530717959E+0000}%
\end{picture}%

%% file: gamman.tex
\unitlength 0.1in
\begin{picture}( 13.0900,  9.8000)( 22.3000,-17.0000)
%
\special{pn 8}%
\special{ar 2844 1366 382 336  0.0000000 6.2831853}%
%
\special{pn 8}%
\special{pa 2506 1212}%
\special{pa 2230 826}%
\special{fp}%
\special{pa 2668 1068}%
\special{pa 2572 720}%
\special{fp}%
\special{pa 3096 1112}%
\special{pa 3318 772}%
\special{fp}%
\special{pa 3212 1284}%
\special{pa 3540 1104}%
\special{fp}%
%
\special{pn 20}%
\special{sh 1}%
\special{ar 2236 834 10 10 0  6.28318530717959E+0000}%
%
\special{pn 20}%
\special{sh 1}%
\special{ar 3216 1280 10 10 0  6.28318530717959E+0000}%
%
\special{pn 20}%
\special{sh 1}%
\special{ar 3092 1116 10 10 0  6.28318530717959E+0000}%
%
\special{pn 20}%
\special{sh 1}%
\special{ar 2672 1068 10 10 0  6.28318530717959E+0000}%
%
\special{pn 20}%
\special{sh 1}%
\special{ar 2506 1208 10 10 0  6.28318530717959E+0000}%
%
\special{pn 20}%
\special{sh 1}%
\special{ar 3530 1104 10 10 0  6.28318530717959E+0000}%
%
\special{pn 20}%
\special{sh 1}%
\special{ar 3314 780 10 10 0  6.28318530717959E+0000}%
%
\special{pn 13}%
\special{sh 1}%
\special{ar 2572 726 10 10 0  6.28318530717959E+0000}%
%
\special{pn 8}%
\special{sh 1}%
\special{ar 2826 922 10 10 0  6.28318530717959E+0000}%
%
\special{pn 8}%
\special{sh 1}%
\special{ar 2750 942 10 10 0  6.28318530717959E+0000}%
%
\special{pn 8}%
\special{sh 1}%
\special{ar 3070 968 10 10 0  6.28318530717959E+0000}%
%
\special{pn 8}%
\special{sh 1}%
\special{ar 2996 934 10 10 0  6.28318530717959E+0000}%
%
\special{pn 8}%
\special{sh 1}%
\special{ar 2912 926 10 10 0  6.28318530717959E+0000}%
\put(27.4100,-13.9700){\makebox(0,0)[lb]{$\Gamma(n)$}}%
\end{picture}%

%% file: gammalambda.tex
\unitlength 0.1in
\begin{picture}( 40.0700,  9.4000)( 19.5700,-37.1000)
%
\special{pn 8}%
\special{pa 4962 3068}%
\special{pa 5374 3068}%
\special{fp}%
%
\special{pn 8}%
\special{pa 4986 3450}%
\special{pa 5380 3450}%
\special{fp}%
%
\special{pn 8}%
\special{sh 0}%
\special{pa 5158 2978}%
\special{pa 5248 2978}%
\special{pa 5248 3574}%
\special{pa 5158 3574}%
\special{pa 5158 2978}%
\special{ip}%
%
\special{pn 8}%
\special{ar 2524 3258 302 318  0.0484118 6.2831853}%
\special{ar 2524 3258 302 318  0.0000000 0.0352158}%
%
\special{pn 20}%
\special{sh 1}%
\special{ar 2764 3448 10 10 0  6.28318530717959E+0000}%
%
\special{pn 20}%
\special{sh 1}%
\special{ar 2756 3078 10 10 0  6.28318530717959E+0000}%
%
\special{pn 20}%
\special{sh 1}%
\special{ar 3716 3708 10 10 0  6.28318530717959E+0000}%
%
\special{pn 20}%
\special{sh 1}%
\special{ar 5252 3068 10 10 0  6.28318530717959E+0000}%
%
\special{pn 20}%
\special{sh 1}%
\special{ar 5924 3688 10 10 0  6.28318530717959E+0000}%
%
\special{pn 20}%
\special{sh 1}%
\special{ar 5172 3068 10 10 0  6.28318530717959E+0000}%
%
\special{pn 20}%
\special{sh 1}%
\special{ar 2508 3578 10 10 0  6.28318530717959E+0000}%
%
\special{pn 20}%
\special{sh 1}%
\special{ar 2508 2948 10 10 0  6.28318530717959E+0000}%
%
\special{pn 20}%
\special{sh 1}%
\special{ar 3148 3078 10 10 0  6.28318530717959E+0000}%
%
\special{pn 20}%
\special{sh 1}%
\special{ar 3100 3258 10 10 0  6.28318530717959E+0000}%
%
\special{pn 20}%
\special{sh 1}%
\special{ar 3156 3448 10 10 0  6.28318530717959E+0000}%
%
\special{pn 20}%
\special{sh 1}%
\special{ar 3572 3528 10 10 0  6.28318530717959E+0000}%
%
\special{pn 20}%
\special{sh 1}%
\special{ar 3700 3258 10 10 0  6.28318530717959E+0000}%
%
\special{pn 20}%
\special{sh 1}%
\special{ar 3588 3018 10 10 0  6.28318530717959E+0000}%
%
\special{pn 20}%
\special{sh 1}%
\special{ar 3724 2778 10 10 0  6.28318530717959E+0000}%
%
\special{pn 20}%
\special{sh 1}%
\special{ar 4260 3270 10 10 0  6.28318530717959E+0000}%
%
\special{pn 20}%
\special{sh 1}%
\special{ar 4436 3268 10 10 0  6.28318530717959E+0000}%
%
\special{pn 20}%
\special{sh 1}%
\special{ar 4732 2958 10 10 0  6.28318530717959E+0000}%
%
\special{pn 20}%
\special{sh 1}%
\special{ar 4732 3578 10 10 0  6.28318530717959E+0000}%
%
\special{pn 20}%
\special{sh 1}%
\special{ar 4972 3448 10 10 0  6.28318530717959E+0000}%
%
\special{pn 20}%
\special{sh 1}%
\special{ar 5164 3448 10 10 0  6.28318530717959E+0000}%
%
\special{pn 20}%
\special{sh 1}%
\special{ar 5244 3448 10 10 0  6.28318530717959E+0000}%
%
\special{pn 20}%
\special{sh 1}%
\special{ar 5396 3448 10 10 0  6.28318530717959E+0000}%
%
\special{pn 20}%
\special{sh 1}%
\special{ar 5780 3538 10 10 0  6.28318530717959E+0000}%
%
\special{pn 20}%
\special{sh 1}%
\special{ar 5924 3268 10 10 0  6.28318530717959E+0000}%
%
\special{pn 20}%
\special{sh 1}%
\special{ar 5788 2998 10 10 0  6.28318530717959E+0000}%
%
\special{pn 20}%
\special{sh 1}%
\special{ar 5932 2778 10 10 0  6.28318530717959E+0000}%
%
\special{pn 20}%
\special{sh 1}%
\special{ar 5324 3268 10 10 0  6.28318530717959E+0000}%
\put(23.0800,-33.5700){\makebox(0,0)[lb]{$\lambda$}}%
%
\special{pn 20}%
\special{ar 2524 3258 262 252  0.0510667 6.2831853}%
\special{ar 2524 3258 262 252  0.0000000 0.0435861}%
\put(59.6400,-36.0700){\makebox(0,0)[lb]{\fbox{$\Gamma'(\lambda)$}}}%
\put(41.1000,-30.5000){\makebox(0,0)[lb]{\fbox{$\Gamma(\lambda)$}}}%
%
\special{pn 8}%
\special{pa 3876 3260}%
\special{pa 4122 3260}%
\special{dt 0.045}%
\special{sh 1}%
\special{pa 4122 3260}%
\special{pa 4056 3240}%
\special{pa 4070 3260}%
\special{pa 4056 3280}%
\special{pa 4122 3260}%
\special{fp}%
%
\special{pn 8}%
\special{pa 5316 3260}%
\special{pa 5920 3260}%
\special{fp}%
%
\special{pn 8}%
\special{pa 3104 3258}%
\special{pa 3708 3258}%
\special{fp}%
%
\special{pn 8}%
\special{pa 5782 3524}%
\special{pa 5920 3686}%
\special{fp}%
%
\special{pn 8}%
\special{pa 5798 2992}%
\special{pa 5942 2770}%
\special{fp}%
%
\special{pn 8}%
\special{pa 4730 2948}%
\special{pa 4730 3568}%
\special{fp}%
%
\special{pn 8}%
\special{pa 4440 3270}%
\special{pa 4244 3270}%
\special{fp}%
%
\special{pn 8}%
\special{ar 5620 3260 302 320  0.0463884 6.2831853}%
\special{ar 5620 3260 302 320  0.0000000 0.0298015}%
%
\special{pn 8}%
\special{ar 4736 3260 302 320  0.0465682 6.2831853}%
\special{ar 4736 3260 302 320  0.0000000 0.0298015}%
\put(20.2000,-30.2700){\makebox(0,0)[lb]{\fbox{$\Gamma$}}}%
%
\special{pn 8}%
\special{pa 3570 3524}%
\special{pa 3708 3682}%
\special{fp}%
%
\special{pn 8}%
\special{pa 3586 3002}%
\special{pa 3730 2780}%
\special{fp}%
%
\special{pn 8}%
\special{pa 2518 2946}%
\special{pa 2518 3566}%
\special{fp}%
%
\special{pn 8}%
\special{pa 2774 3450}%
\special{pa 3168 3450}%
\special{fp}%
%
\special{pn 8}%
\special{pa 2758 3068}%
\special{pa 3170 3068}%
\special{fp}%
%
\special{pn 8}%
\special{ar 3408 3258 302 318  0.0483183 6.2831853}%
\special{ar 3408 3258 302 318  0.0000000 0.0352158}%
%
\special{pn 8}%
\special{pa 2220 3260}%
\special{pa 1974 3260}%
\special{fp}%
%
\special{pn 20}%
\special{sh 1}%
\special{ar 1960 3260 10 10 0  6.28318530717959E+0000}%
%
\special{pn 20}%
\special{sh 1}%
\special{ar 2220 3260 10 10 0  6.28318530717959E+0000}%
%
\special{pn 8}%
\special{sh 0}%
\special{pa 4650 3180}%
\special{pa 4810 3180}%
\special{pa 4810 3370}%
\special{pa 4650 3370}%
\special{pa 4650 3180}%
\special{ip}%
%
\special{pn 20}%
\special{sh 1}%
\special{ar 4730 3180 10 10 0  6.28318530717959E+0000}%
%
\special{pn 20}%
\special{sh 1}%
\special{ar 4730 3370 10 10 0  6.28318530717959E+0000}%
\end{picture}%